\documentclass[12pt]{article}

\font\goth = eufm10 scaled \magstep 1
 2
\font\Bbb = msbm10 scaled \magstep 1
\def\UU{{\cal U}}
\def\MM{{\cal M}}
\def\FF{{\cal F}}
\def\RR{{\cal R}}
\def\gtg{\mbox{\goth g}}
\def\BN{\mbox{\Bbb N}}
\def\BC{\mbox{\Bbb C}}
\def\BR{\mbox{\Bbb R}}
\def\BR{\mbox{\Bbb R}}
\def\veps{\varepsilon}
\def\nat{\natural}
\def\diag{\mbox{diag}}
\def\sgn{\mbox{sgn}}
\def\span{\mbox{span}}
\def\nonum{\nonumber}
\def\slr{\mbox{\goth sl}(2,\BR)}
\def\slc{\mbox{\goth sl}(2,\BC)}
\def\exp{\mathrm{exp}}
\def\I{\mathrm{i}}

\begin{document}
\date{}
\author{\textsc{Pavel \v S\v tov\'\i\v cek}\\ 
\makebox{}\\
 {\itshape Department of Mathematics}\\
 {\itshape Faculty of Nuclear Science, Czech Technical University}\\
 {\itshape Trojanova 13, 120 00 Prague, Czech Republic}\\
 {\itshape stovicek@km1.fjfi.cvut.cz}
}
\title{Some remarks on $\UU_{q}({\boldmath sl}(2,{\boldmath R}))$ 
at root of unity}
\maketitle

\begin{abstract}
We discuss a modification of $\UU_q(\slr)$ and a class of its 
irreducible representations when $q$ is a root of unity.
\end{abstract}

\section{Introduction}

Nowadays \( q \)-deformed universal enveloping algebras 
\( \UU _{q}(\gtg ) \)
are understood in depth in the case when \( \gtg  \) is a complex simple Lie
algebra belonging to one of the four principal series. The same is true for
compact forms of these Lie algebras (see, e.g., monographs
\cite{[1]}, \cite{[2]}, \cite{[3]}). 
On the other hand, attempts to introduce \( q \)-deformed enveloping
algebras for non-compact real Lie algebras frequently lead to 
serious difficulties 
though several particular cases have been already studied (see, e.g.,
\cite{[4]}, \cite{[5]}, \cite{[6]}). 
In this note we discuss one of the simplest examples with
\( \gtg =\slr  \) as a real form of \( \slc  \). The deformation parameter
\( q \) is supposed to be a root of unity,
\[
q=\exp (\I \pi P/Q),
\]
where \( Q\in \BN  \) is odd, \( P\in \{1,\ldots ,Q-1\} \), and \( P \) and
\( Q \) are relatively prime integers. So $q^{2j}\neq1$, 
$j=1,\ldots,Q-1$, and $q^{2Q}=1$.

We use the standard definition of the Hopf algebra \( \UU _{q}(\slc ) \) 
with the generators \( K,\, K^{-1},\, E,\, F \), the defining relations
\begin{eqnarray*}
  & K\, K^{-1}=K^{-1}K=1,\textrm{ }K\, E=q\, E\, K,
  \textrm{ }K\, F=q^{-1}F\, K, & \\
  & [\, E,F\, ]=\frac{1}{q-q^{-1}}(K^{2}-K^{-2}), & 
\end{eqnarray*}
the comultiplication
\[
\Delta K=K\otimes K,\textrm{ }\Delta E=K\otimes E+E\otimes K^{-1},
\textrm{ }\Delta F=K\otimes F+F\otimes K^{-1},
\]
the antipode
\[
S(K)=K^{-1},\textrm{ }S(E)=-q^{-1}E,\textrm{ }S(F)=-q\, F,
\]
and the counit
\[
\varepsilon (K)=1,\textrm{ }\varepsilon (E)=\varepsilon (F)=0.
\]

A real form is determined by a \( \ast  \)-involution; 
an element \( X \) of a complex Hopf algebra 
belongs to a real form if and only if \( X^{\ast }=S(X) \). 
Particularly, \( \UU _{q}(\slr ) \)
is determined by the \( \ast  \)-involution
\begin{equation}
  \label{star_invol}
  K^{\ast }=K,\textrm{ }E^{\ast }=-q^{-1}E,\textrm{ }F^{\ast }=-q\, F.
\end{equation}
Necessarily, \( q \) is a complex unit, \( \bar{q}=q^{-1} \). 

Usually it is more convenient to deal with the complexification 
of a real form. In that case one
regards the real form as the original complex Hopf algebra but
endowed, in addition, 
with the \( \ast  \)-involution in question. 
We shall adopt this point of view and treat \( \UU _{q}(\slr ) \) 
as the complex Hopf algebra \( \UU _{q}(\slc ) \) with the 
$\ast$-involution (\ref{star_invol}).

\section{A modification of \protect\( \UU _{q}(\slr )\protect \)}

Let \( \UU  \) be a \( \ast  \)-Hopf subalgebra of 
\( \UU _{q}(\slr ) \) generated
by \( X,\, Y,\, Z,\, Z^{-1} \), where
\[
X=-\I \, q^{-1}E\, K^{-1},\textrm{ }Y=-\I \, q\, F\, K^{-1},
\textrm{ }Z=K^{-2}.
\]
Thus \( \UU  \) is defined by the relations
\begin{equation}
\label{defrel_Unat}
Z\, X=q^{-2}X\, Z,\textrm{ }Z\, Y=q^{2}Y\, Z,
\textrm{ }q^{-1}X\, Y-q\, Y\, X=-\frac{1}{q-q^{-1}}(1-Z^{2}),
\end{equation}
with the comultiplication
\[
\Delta Z=Z\otimes Z,\textrm{ }\Delta X=1\otimes X+X\otimes Z,
\textrm{ }\Delta Y=1\otimes Y+Y\otimes Z,
\]
the antipode
\[
S(Z)=Z^{-1},\textrm{ }S(X)=-X\, Z^{-1},\textrm{ }S(Y)=-Y\, Z^{-1},
\]
and the counit
\[
\veps (Z)=1,\textrm{ }\veps (X)=\veps (Y)=0.
\]
Furthermore, all the generators are Hermitian,
\[
Z^{\ast }=Z,\textrm{ }X^{\ast }=X,\textrm{ }Y^{\ast }=Y.
\]
It is also straightforward to check that
$$
C=X\,Y\,Z^{-1}-\frac{1}{(q-q^{-1})^2}(Z^{-1}+q^2Z)
$$
is an Casimir element in $\UU$.

Unfortunately, there exists no non-trivial irreducible representations
\( \rho  \) of \( \UU  \). 
Actually, \( Z^{Q} \) belongs to the center of \( \UU  \)
and is Hermitian. Thus, by the Schur lemma, \( \rho (Z)^{Q}=c\, I \) for some
real \( c\neq 0 \). Consequently, the self-adjoint operator \( \rho (Z) \) 
is a multiple of the identity as well. The commutation relations then 
imply that \( \rho (X)=\rho (Y)=0 \), \( \rho (Z)=\pm I \).

To improve this situation we propose a modification of \( \UU  \) that we call
here \( \UU ^{\natural } \). As a Hopf algebra, \( \UU  \) is extended to
\( \UU ^{\natural } \) by adding another generator, \( T \), which satisfies
\[
T^{2}=1,\textrm{ }\Delta T=T\otimes T,\textrm{ }S(T)=T,
\textrm{ }\veps (T)=1.
\]
A \( \ast  \)-involution on \( \UU ^{\natural } \) is defined as follows:
\[
X^{\ast }=T\, X\, T,\textrm{ }Y^{\ast }=
T\, Y\, T,\textrm{ }Z^{\ast }=T\, Z\, T,\textrm{ }T^{\ast }=T.
\]
So \( \UU  \) is a Hopf subalgebra of \( \UU ^{\natural } \) but not a
\( \ast  \)-Hopf subalgebra. On the other hand, 
\( \UU  \) may be obtained from \( \UU ^{\natural } \)
by specializing \( T \) to 1.

\section{A class of representations of \protect\( \UU ^{\nat }\protect \)}

Next we present a class of irreducible representations of the 
\( \ast  \)-algebra \( \UU ^{\nat } \) while the question of 
a complete classification of irreducible
representations of \( \UU ^{\nat } \) is proposed as an open problem. Though
it is not excluded that the definition of \( \UU ^{\nat } \) should be further
modified in order to get a reasonable theory. 
In this section most steps are only outlined with some details omitted.

The representation \( \rho  \) depends on an integer parameter 
\( n\in \{1,2,\dots ,Q\} \)
and its dimension \( d \) equals \( Q+1-n \). The matrices 
\( \rho (X),\, \rho (Y),\, \rho (Z) \)
are tridiagonal with non-vanishing entries\-
\begin{eqnarray}
  \rho (Z)_{m-1,m} & = & -(q^{2}-q^{-2})q^{2m+n-1}a_{m},\nonum \\
  \rho (Z)_{m+1,m} & = & (q^{2}-q^{-2})q^{-2m-n-1}b_{m+1},\label{rhoZ} \\
  \rho (Z)_{mm} & = & (q+q^{-1})c_{m},\nonum 
\end{eqnarray}
\begin{eqnarray}
  \rho (X)_{m-1,m} & = & (q+q^{-1})a_{m},\nonum \\
  \rho (X)_{m+1,m} & = & (q+q^{-1})b_{m+1},\label{rhoX} \\
  \rho (X)_{mm} & = & d_{m},\nonum 
\end{eqnarray}
\begin{eqnarray}
  \rho (Y)_{m-1,m} & = & (q+q^{-1})q^{2(2m+n-1)}a_{m},\nonum \\
  \rho (Y)_{m+1,m} & = & (q+q^{-1})q^{-2(2m+n+1)}b_{m+1},\label{rhoY} \\
  \rho (Y)_{mm} & = & -d_{m},\nonum 
\end{eqnarray}
\( m=0,1,\dots \, d-1 \). Here
\begin{eqnarray*}
  a_{m} & = & b_{m}\\
  & = & \frac{1}{q^{2m+n-1}+q^{-2m-n+1}}
  \sqrt{\frac{[m]_{q^{2}}[m+n-1]_{q^{2}}}{(q^{2m+n-2}+q^{-2m-n+2})
      (q^{2m+n}+q^{-2m-n})}},\\
  c_{m} & = & \frac{q^{n-1}+q^{-n+1}}{(q^{2m+n-1}+q^{-2m-n+1})
    (q^{2m+n+1}+q^{-2m-n-1})},\\
  d_{m} & = & \frac{q^{n-1}+q^{-n+1}}{(q^{2m+n-1}+q^{-2m-n+1})
    (q^{2m+n+1}+q^{-2m-n-1})}[2m+n]_{q}.
\end{eqnarray*}
The quantum numbers are defined as usual, 
\[ 
[x]_{q}=\frac{q^{x}-q^{-x}}{q-q^{-1}}.
\]
The matrix \( \rho (T) \) is diagonal,
\[
\rho (Z)=\diag (\tau _{0},\tau _{1},\dots ,\tau _{d-1})\]
where \( \tau _{0}=1 \) and
\begin{equation}
  \label{sign}
  \frac{\tau _{m}}{\tau _{m-1}}=\sgn \left( 
    \frac{[m]_{q^{2}}[m+n-1]_{q^{2}}}{(q^{2m+n-2}+q^{-2m-n+2})
(q^{2m+n}+q^{-2m-n})}\right) ,
\end{equation}
\( m=1,2,\dots ,d-1 \).

Let us remark that a source of difficulties when working with real 
forms comes from the fact that the deformation parameter $q$ is forced
to be a complex unit. In that case the sign $\tau_m/\tau_{m-1}$ in
(\ref{sign}) may equal -1 for particular values of $m$. Concerning the
representation $\rho$, it is worth mentioning that the matrix 
$\rho((q\,X-q^{-1}Y)Z^{-1})$ is diagonal and
$$
\rho((q\,X-q^{-1}Y)Z^{-1})_{mm}=[2m+n]_q.
$$

The verification of the commutation relations (\ref{defrel_Unat}) 
is straightforward. This may be done even in the case when $q$ is
generic and the tridiagonal matrices (\ref{rhoZ}), (\ref{rhoX}), 
(\ref{rhoY}) are infinite with $m=0,1,2,\ldots$. 
One then finds that relations (\ref{defrel_Unat}) 
are satisfied if and only if the the coefficients
\( c_{m} \) obey a recursive equation,
\begin{eqnarray*}
  &  & (q^{2m+n+3}+q^{-2m-n-3})c_{m+1}-(q+q^{-1})(q^{2m+n}+q^{-2m-n})c_{m}\\
  &  & +(q^{2m+n-3}+q^{-2m-n+3})c_{m+1}=0,
\end{eqnarray*}
and \( d_{m} \) and \( a_{m}b_{m} \) are expressed in terms of \( c_{m} \),
\begin{eqnarray*}
  d_{m} & = & -\frac{(q^{2m+n-1}+q^{-2m-n+1})c_{m}-
    (q^{2m+n-3}+q^{-2m-n+3})c_{m-1}}{(q-q^{-1})^{2}},\\
  a_{m}b_{m} & = & (q-q^{-1})^{-4}(q+q^{-1})^{-2}\\
  &  & \times (q^{2m+n}+q^{-2m-n})^{-1}(q^{2m+n-2}+q^{-2m-n+2})^{-1}\\
  &  & \times \big((q^{2m+n+1}+q^{-2m-n-1})^{2}c^{\, 2}_{m}+
  (q^{2m+n-3}+q^{-2m-n+3})^{2}c^{\, 2}_{m-1}\\
  &  & \, -(q^{2}+q^{-2})(q^{2m+n-3}+q^{-2m-n+3})(q^{2m+n+1}+
  q^{-2m-n-1})c_{m}c_{m-1}\\
  &  & \, +(q-q^{-1})^{2}\big ).
\end{eqnarray*}
Equivalently,
\[
d_{m}=-\frac{(q^{2m+n+3}+q^{-2m-n-3})c_{m+1}-(q^{2m+n+1}+q^{-2m-n-1})c_{m}}
{(q-q^{-1})^{2}}.
\]

To verify the irreducibility we shall show that even the restriction
of \( \rho  \)
to the subalgebra \( \UU  \) is irreducible. This will become obvious as soon
as we prove that \( \rho  \) is equivalent to \( \tilde{\rho } \) with
\begin{eqnarray*}
  & \tilde{\rho }(X)=\left( \begin{array}{ccccc}
      0 & 0 & \ldots  & 0 & 0\\
      x_{1} & 0 & \ldots  & 0 & 0\\
      0 & x_{2} & \ldots  & 0 & 0\\
      \vdots  & \vdots  & \ddots  & \vdots  & \vdots \\
      0 & 0 & \ldots  & x_{d-1} & 0
    \end{array}\right) ,
  \textrm{ }\tilde{\rho }(Y)=\left( 
    \begin{array}{ccccc}
      0 & y_{1} & 0 & \ldots  & 0\\
      0 & 0 & y_{2} & \ldots  & 0\\
      \vdots  & \vdots  & \vdots  & \ddots  & \vdots \\
      0 & 0 & 0 & \ldots  & y_{d-1}\\
      0 & 0 & 0 & \ldots  & 0
    \end{array}\right) , & 
\end{eqnarray*}
and
\begin{eqnarray*}
  & \tilde{\rho }(Z)=\left( \begin{array}{cccc}
      z_{0} & 0 & \ldots  & 0\\
      0 & z_{1} & \ldots  & 0\\
      \vdots  & \vdots  & \ddots  & \vdots \\
      0 & 0 & \ldots  & z_{d-1}
    \end{array}\right) ,  & 
\end{eqnarray*}
where
\[
x_{j}=q^{-n-2j+1}[n+j-1]_{q},\textrm{ }y_{j}=[j]_{q},
\textrm{ }z_{j}=q^{-n-2j}.
\]
Note that $x_j\neq0$, $y_j\neq0$, for $j=1,\ldots,d-1$.

The equivalence in turn follows from a more geometrical realization 
of the representation \( \rho  \) which is closely related to 
the twisted adjoint action \cite{[7]}, \cite{[8]}. 
The vector space \( \MM  \) of meromorphic functions
in variable \( w \) on the complex plane becomes a left \( \UU  \) module
with respect to the action
\begin{eqnarray*}
  X\cdot f(w) & = & 
  -\I\frac{q^{-1}w}{q-q^{-1}}
  \left( q^{n}f(w)-q^{-n}f(q^{-2}w)\right) ,\\
  Y\cdot f(w) & = & 
  \I\frac{q^{-n+1}}{(q-q^{-1})w}\left( f(w)-f(q^{-2}w)\right) ,\\
  Z\cdot f(w) & = & q^{-n}f(q^{-2}w).
\end{eqnarray*}
Set
\[
\psi _{m}(w)=\frac{\prod ^{m-1}_{j=0}(q^{2j+n}w-\I )}
{\prod ^{n+m-1}_{j=0}(q^{2j-2m-n}w+\I )},\quad m=0,1,2,\dots .
\]
Then the vector space 
\( \span \{\psi _{0},\psi _{1},\ldots ,\psi _{d-1}\} \),
\( d=Q+1-n \), may be checked to be \( \UU  \) invariant. 
After renormalization of the basis
vectors, \( \tilde{\psi }_{m}=\lambda _{m}\psi _{m} \), 
with the factors \( \lambda _{m} \)
being determined by \( \lambda _{0}=1 \) and
\[
\lambda _{m}=\sqrt{
\frac{(q^{2m+n}+q^{-2m-n})[m+n-1]_{q^{2}}}
{(q^{2m+n-2}+q^{-2m-n+2})[m]_{q^{2}}}\, }\lambda _{m-1},
\]
we get the representation \( \rho  \). 

Consider now a point set, 
\[ 
M=\{1,q^{2},\ldots ,q^{2Q-2}\}\subset \BC .
\]
Note that for any function \( f \), the values of the function \( A\cdot f \)
on the set \( M \) depend only on the restriction \( f|_{M} \) where \( A \)
is any of the generators \( X \), \( Y \) or \( Z \). Thus the vector space
\( \FF \cong \BC ^{Q} \) of functions on \( M \) becomes a \( \UU  \) module
and the restriction map 
\( \MM\to\FF:f\mapsto f|_{M} \) 
is a surjective morphism of  \( \UU  \) modules. 
The representation \( \rho  \) corresponds to the submodule 
\( \RR =\span \{\tilde{\psi }_{0}|_{M},
\tilde{\psi }_{1}|_{M},\ldots ,\tilde{\psi }_{d-1}|_{M}\} \)
with the distinguished basis.
Omitting the details we claim that another basis in 
\( \RR  \) may be chosen as 
\( \{\phi _{0}|_{M},\phi _{1}|_{M},\ldots ,\phi _{d-1}|_{M}\} \)
where
\[
\phi _{j}(w)=
\left( \frac{1}{\I }\right) ^{j}q^{\frac{1}{2}j(j-1)+nj}w^{-Q+j}.
\]
Expressing operators in the latter basis we get the representation 
\( \tilde{\rho } \).
This proves the equivalence of \( \rho  \) and \( \tilde{\rho } \) 
and consequently that the representation $\rho$ is irreducible.


\end{document}